\documentclass [11pt]{amsart} 
\usepackage {amsmath, amssymb, amscd, amsthm, epsfig, color}

\input xy
\xyoption {all}

\def \bb {\mathbb}

\def \n {{{\mathcal {N}}}}
\def \m {{{\mathcal {M}}}}

\def \e {{{\mathcal {E}}}}
\def \f {{{\mathcal {F}}}}

\def \k {{{\mathcal {K}}}}

\def \t {{{\mathcal {T}}}}
\def \mor {{{\text{Mor}}_{d}(C, G(r,N))}}
\def \quot {{{\text{Quot}}_{d} ({\mathcal O}^N, r, C)}}
\def \gbar {{\bar g}}
\def \p {{\mathbb P}}

\newtheorem{theorem}{Theorem} 
\newtheorem {lemma}{Lemma} 
 
\newtheorem {proposition}{Proposition}

\theoremstyle{definition}
\newtheorem{remark}{Remark}

\theoremstyle {definition}

\begin{document}

\title[Virtual intersections on the Quot-scheme]{Virtual intersections on the Quot scheme and Vafa-Intriligator formulas}
\author {Alina Marian}
\address {Department of Mathematics}
\address {Yale University}
\email {alina.marian@yale.edu}
\author {Dragos Oprea}
\address {Department of Mathematics}
\address {Massachusetts Institute of Technology.}
\email {oprea@math.mit.edu}
\date{} 

\begin {abstract}

We construct a virtual fundamental class on the Quot scheme parametrizing quotients of a trivial bundle on a smooth projective curve. We use the virtual localization formula to calculate virtual intersection numbers on Quot. As a consequence, we reprove the Vafa-Intriligator formula; our answer is valid even when the Quot scheme is badly behaved. More intersections of Vafa-Intriligator type are computed by the same method. Finally, we present an application to the non-vanishing of the Pontrjagin ring of the moduli space of bundles. 

\end {abstract}

\maketitle

\section{Introduction}

We study the intersection theory of the Quot scheme $\quot$ of degree $d$, rank $N-r$ coherent sheaf quotients of ${\mathcal {O}}^N$ on a smooth complex projective curve $C$ of genus $g$ via equivariant localization.

$\quot$ provides a compactification of the scheme $\mor$ of degree $d$ morphisms from $C$ to the Grassmannian $G(r, N),$ and has been analyzed a lot
from this point of view. Indeed, $\mor$ sits inside $\quot$ as the subscheme of {\em{locally free}} quotients of ${\mathcal {O}}^N$. For large
degree $d$, certain counts of maps from $C$ to $G(r,N)$ can in fact be realized as intersection numbers on $\quot$, and can be carried out in this
setting \cite{bertram1}, \cite{bertram}. Part of the intersection theory on $\quot$ has thus been well studied. In this paper, we note that all of the
intersection theory on $\quot$ can be explicitly computed by exploiting the natural ${{\bb C}^{\star}}$-action of the situation and the ubiquitous
virtual localization theorem \cite{graberpandharip}. In particular we recover the intersection numbers which have already been known, and which are
given by the elegant formula of Vafa and Intriligator.

We now detail the discussion. Quot schemes have been shown by Grothendieck \cite{gro}, in 
all
generality, to be fine moduli spaces for the problem of parametrizing quotients of a fixed 
sheaf, and as such to carry universal structures. For the case under study, quotients of 
${\mathcal {O}}^N$ on a curve $C$, let $$0 \rightarrow \e \rightarrow {\mathcal {O}}^N \rightarrow \f \rightarrow 0$$ be the
universal family on $\quot \times C.$ We note for further use that the Zariski tangent 
space to $\quot$ at a closed point $[0\rightarrow E \rightarrow {\mathcal {O}}^N \rightarrow F \rightarrow 0]$ 
is
$\text {Ext}^0(E, F)$, while the obstructions lie in $\text {Ext}^1(E, F)$. Hence 
the
expected dimension of $\quot$ is $$e = \chi (E^{\vee} \otimes F) = Nd -r(N-r)\gbar.$$ As usual, we write $\gbar=g-1$.

Now let $\{1, \delta_k, 1\leq k \leq 2g, \omega \}$ be a symplectic basis for the 
cohomology of $C$ and
let $$c_i(\e^{\vee}) = a_i\otimes 1 + \sum_{j=1}^{2g} b_i^j \otimes \delta_j  + f_i \otimes 
\omega$$ be
the K\"{u}nneth decomposition of the Chern classes of the dual universal bundle $\e^{\vee}.$ Note that $$a_i
\in H^{2i}(\quot),\;\; b_i^j \in H^{2i-1}(\quot), f_i \in H^{2i-2}(\quot).$$ Moreover, let $p$ be a 
point
on $C$ and let $\e_p$ denote the restriction of $\e$ to $\quot \times \{p\}$. We write $$\eta: \quot \times 
C
\rightarrow \quot$$ for the projection on the first factor. It is clear that $$a_i = c_i (\e_{p}^{\vee}) \, \, \, \text{and}
\, \, f_i = \eta_{\star} c_i(\e^{\vee}).$$

Since the Quot scheme may have several possibly oversized irreducible components, making sense of intersection numbers can be subtle. The machinery of virtual
moduli cycles of \cite {litian}, \cite {behrendfantechi} was developed precisely to deal with such issues. Our first result is the existence of a virtual fundamental class on $\quot$. It will sometimes be convenient to write $Q_d$ for $\quot$, when no confusion is likely. 

\begin {theorem} The scheme $\quot$ admits a perfect obstruction theory and a virtual fundamental class $[Q_d]^{vir}$ of the expected dimension
$e$. \end {theorem}

This should come as no surprise, and has in fact been known to some experts. The referee pointed out that an independent 
proof of this result is contained in the unpublished manuscript \cite{CFK2}. There, the construction of the virtual fundamental 
class relies on the existence, demonstrated in \cite{CFK1}, of a dg-manifold, the derived Quot scheme, whose degree $0$ truncation is the 
usual Quot scheme. Moreover, constructions of virtual fundamental classes of similar flavor arise in higher dimensions in the context of Seiberg-Witten theory for the Hilbert scheme of divisors on a surface \cite {okonek} and Donaldson-Thomas theory 
for moduli spaces of sheaves on threefolds \cite {thomas}. We remark that our arguments can be extended to more general 
situations; in particular it is not necessary to consider quotients of the trivial bundle. However, this case suffices for the 
applications we have in mind.

Next, fixing $p\in C$, there is a natural embedding $$i_p: \quot\longrightarrow \text 
{Quot}_{d+r}({\mathcal {O}}^N,
r, C)$$ obtained by twisting subbundles by ${\mathcal {O}}(-p)$. This exhibits $\quot$ as the zero locus of a section of $\left(\mathcal E_p^{\oplus
N}\right)^{\vee}$. The virtual classes we construct are compatible with respect to this embedding.

\begin {theorem} \label{compatible} The following equality holds true:
$$i_{\star}\left[\quot\right]^{vir}=c_{\text {top}} (\mathcal E^{\vee}_p)^N \cap \left[\text
{Quot}_{d+r}({\mathcal {O}}^N, r, C)\right]^{vir}.$$ \end {theorem}

Now $GL(N)$ acts naturally on $\quot$ by sending quotients on $C$ to their
composition with automorphisms of ${\mathcal {O}}^N$. Furthermore, the dual universal bundle $\e^{\vee}$ is
$GL(N)$-equivariant. We will look at the restriction of the action to a diagonal ${\bb C}^{\star}$
inside $GL(N),$ and will apply the virtual localization formula of \cite{graberpandharip}, obtaining that

\begin {proposition}\label{allint} All monomials in the $a$, $b$ and $f$ classes can be explicitly
evaluated on the virtual fundamental class $[Q_d]^{vir}$. \end {proposition}

For large degree $d$ relative to $r$, $N$ and $g$, $\quot$ is generically reduced, irreducible
of the expected dimension \cite{bertram1}. The virtual fundamental class agrees with the usual fundamental class,
hence we are computing honest intersection numbers. In this large-degree regime, intersections of 
$a$-{\em{classes}} were studied extensively in \cite{bertram}. It is shown there that 
the evaluation of a top-degree monomial in the $a$-classes on the fundamental cycle has enumerative meaning: it is the number of degree $d$ maps from $C$ to the Grassmannian $G(r,N)$ which send fixed points $p_i$ on $C$ to special Schubert
subvarieties of $G(r,N).$ The number of such maps is computed by the Vafa-Intriligator formula 
\cite{intrilig}. 
Further in \cite{bertram}, the author evaluates, up to a calculation in genus 0, the top intersection of 
$a$-classes on $\quot$ by degenerating to lower genus, hence he essentially proves the Vafa-Intriligator
formula. In addition, he defines $a$-intersection numbers in low degree. Our Theorem \ref{compatible}
shows that Bertram's definition gives rise to virtual intersection numbers. Thus our approach 
clarifies issues 
related to the bad behavior of the moduli space. 

It should also be mentioned that a complete proof of the Vafa-Intriligator formula is contained in \cite{sieberttian}; it exploits the standard presentation of the small quantum cohomology ring of $G(r,N)$ and uses too a degeneration of the curve to genus 0.

Using localization, we explicitly compute the $a$-intersection numbers, recovering the formula  
of Vafa and Intriligator. The formulation below follows Bertram's paper \cite {bertram}. In the proof, we explain how
each term in the Vafa-Intriligator formula arises as a contribution from various fixed loci. We will
investigate a similar formula for HyperQuot schemes in future work.

\begin {theorem} \label{vafaint}Let $P(X_1, \ldots, X_r)$ be a polynomial of weighted degree $e$, where
the variable $X_i$ has degree $i$. Define $$J(x_1, \ldots, x_r)=N^r \cdot 
x_1^{-1}\cdot \ldots \cdot
x_r^{-1}\cdot \prod_{1\leq i<j\leq r}(x_i-x_j)^{-2}.$$ Then, $$P(a_1, \ldots, a_r)\cap
\left[Q_d\right]^{vir}=u\cdot\sum_{\lambda_1, \ldots, \lambda_r} 
R(\lambda_1,
\ldots, \lambda_r) J^{g-1}(\lambda_1, \ldots, \lambda_r),$$ the sum being taken over all $\binom{N}{r}$
tuples $(\lambda_1, \ldots, \lambda_r)$ of distinct $N^{\text th}$ roots of unity. Here, $$u=(-1)^{(g-1)\binom{r}{2}+d(r-1)},$$ and $R$ is the symmetric polynomial
obtained by expressing $P(a_1, \ldots, a_r)$ in terms of the Chern roots of $\e^{\vee}_p$. 
\end {theorem}

Moreover, we will strengthen the theorem above by considering intersections of $a$ and certain $b$ classes. While more general formulas can be written down, we will only note here that
\begin {proposition}\label{bees}
For all $s\leq d$ and $1\leq j_1 <\ldots <j_s\leq g$,  
$$\left(b_1^{j_1}b_1^{j_1+g} \cdot \ldots\cdot b_1^{j_s} b_1^{j_s+g}\right)\, P(a)\cap \left[Q_d\right]^{vir}=\frac{u} {N^{s}} \sum_{\lambda_1, \ldots, \lambda_r} (\lambda_1+\ldots+\lambda_r)^s \cdot \left(RJ^{\gbar}\right)(\lambda_1, \ldots, \lambda_r).$$ The expression evaluates to $0$ for $s>d$.
\end {proposition}

We observed in Proposition $\ref{allint}$ that all intersection numbers on $Q_d$ can in
principle be computed; precise formulas are not always easy to obtain. We further derive the following
intersection number which is reminiscent of the expression obtained by Vafa and Intriligator as well.

\begin {theorem} \label {vinew} Let $P$ be any polynomial of weighted degree $e-l+1$, where $2\leq l\leq r$. Let $R$ be the symmetric polynomial obtained by expressing $P(a_1, \ldots, a_r)$ in terms of the Chern roots. The intersection product $$f_l\cdot P(a_1, \ldots, a_r)\cap \left[Q_d\right]^{vir}$$
equals the following sum over all the $\binom{N}{r}$ tuples of distinct $N$-roots of unity
$$\frac{u}{N}\sum_{\lambda_1, \ldots, \lambda_r}(\mathcal D_l
R)(\lambda_1, \ldots, \lambda_r) \cdot J^{g-1}(\lambda_1, \ldots, \lambda_r),$$ where $\mathcal D_l$ is
a first order linear differential operator defined by equation $\eqref{tl}.$ \end {theorem}

An intersection number on $\quot$ involving $f$ classes
should morally correspond to a count of degree $d$ maps from $C$ to $G(r,N)$ whose images intersect 
given
special Schubert subvarieties of $G(r,N).$ Unfortunately, due to the bad nature of the boundary of the  
$\quot$ compactification, an $a$ and $f$ intersection number does not have enumerative
meaning in general; some of these issues are addressed in \cite {marian}.

Nonetheless, evaluating general top monomials in $a$, $f$ and $b$ classes has important applications
to the computation of the intersection numbers on the moduli space $\mathcal M(r, d)$ of rank $r$,
degree $d$ bundles on $C$, when $r$, $d$ are coprime. It was shown in \cite{marian} that for large $N$
and $d$, any top intersection product of the cohomology generators of $\mathcal M (r,d)$ can be 
realized
as a top intersection of $a$, $f$ and $b$ classes on $\quot$, hence is computable by the methods of
this work. This gives an algebraic way of calculating the intersection theory of $\mathcal M(r,d)$, 
and could thus provide an alternative to the pioneering work \cite{jeffreykirwan}. An ampler 
exploration of 
the intersection theory on $\mathcal M(r,d)$ via equivariant localization on $\quot$ will be taken up in a future paper. Note that complete results in rank $2$ are obtained in \cite {opreamarian}.

The vanishing of the ring generated by the Pontrjagin classes of $\mathcal M(r,d)$ in degrees higher than $r(r-1)(g-1)$ was proved in \cite{kirwan}. There, the optimality of the result is also shown, by producing a nonzero element of degree $r(r-1)(g-1)$. We will replicate this optimality statement as an application of the  
Vafa-Intriligator formula alone. 

\begin {theorem}\label{modbun} There exists a non-zero element of degree $r(r-1)(g-1)$ in the 
Pontrjagin ring
of $\mathcal M(r, d)$. \end {theorem}

The paper proceeds as follows. In Section 2, we show that $\quot$ standardly admits a perfect obstruction
theory compatible with the natural ${\bb C}^{\star}$ action, hence that it has a virtual fundamental class
$[\quot]^{vir} \in H_{2e}(\quot)$ compatible with the action. In Section 3, we describe the fixed loci of
the torus action and the weights on the virtual normal bundles to the fixed loci. This analysis was carried out
already in \cite{strom} for the case of the Quot scheme on ${\bb P}^1$ and is no different in the case
when $C$ has arbitrary genus. We further calculate the inverses of the equivariant Euler classes of the
virtual normal bundles, and prove the compatibility statement expressed by Theorem $\ref {compatible}$. In the
fourth section we obtain explicit formulas for the virtual intersections of $a$ classes, recovering the
formula of Vafa and Intriligator. In the fifth section, we indicate how we explicitly evaluate any top
degree monomial and prove Theorem $\ref{vinew}$. Finally, we prove Theorem $\ref{modbun}$ in the last
section.\\

{\bf Acknowledgements.} We would like to thank Martin Olsson, Rahul Pandharipande, Johan de Jong, and Max Lieblich 
for conversations related to Section $2$ of this work. Moreover, we thank the referee for alerting us to the 
unpublished preprint \cite{CFK2}, and Mikhail Kapranov for making this manuscript readily available to us.

\section {The perfect obstruction theory of the Quot scheme}

We set out to note the existence of a virtual fundamental class of the expected dimension $e$ on the Quot scheme $Q_d$.

\addtocounter{theorem}{-5} \begin {theorem} The Quot scheme of quotients of the trivial bundle on
a smooth projective curve admits a perfect obstruction theory and a virtual fundamental class of the expected
dimension $e$. \end {theorem}

{\it Proof.} There are essentially two methods of constructing virtual fundamental classes: they are due to
Li-Tian \cite {litian} and Behrend-Fantechi \cite {behrendfantechi} respectively. The
Li-Tian approach is more naturally suited to the problem at hand; the discussion below will use this
approach as the starting point.

To construct the virtual fundamental class, it suffices to produce a two-step complex of vector
bundles $\left[\mathcal A_0\to \mathcal A_1\right]$ on $Q_d$ resolving the tangent-obstruction
complex. That is, for all sheaves $\mathcal I$ on $Q_d$, we would like an exact sequence
\begin{equation}\label{tanobs} 0\to \mathcal T^{1} (\mathcal I)\to \mathcal
A_0\otimes \mathcal I\to \mathcal A_1\otimes \mathcal I \to \mathcal T^2 \otimes \mathcal I\to 0.\end {equation} Here, we write $\mathcal T^1 (\mathcal I)$ and $\mathcal T^2\otimes \mathcal I$ for the deformation and obstruction spaces on $Q_d$ in the direction of $\mathcal I$. The precise definitions are contained in \cite {litian}. In the case of the Quot scheme, it is well known that \cite {flenner} $$\mathcal T^1(\mathcal I)=Hom_{\eta}(\mathcal E, \mathcal F \otimes \eta^{\star}\mathcal I),\;\; \mathcal T^2 = Ext^1_{\eta}(\mathcal E, \mathcal F).$$

To produce the resolution $\eqref{tanobs}$, we recall the Grothendieck embedding of the Quot
scheme into a Grassmannian.  Let $m$ be a positive integer, large enough so that $\e(m)$ and $\f
(m)$ are fiberwise generated by sections and $$R^1 \eta_{\star}\e(m)=R^1\eta_{\star} \f(m)=0.$$
Then $R^{0}\eta_{\star}\e (m)$ and $R^{0}\eta_{\star}\f (m)$ are locally free. Note that here we
use ${\mathcal {O}}(m)$ to denote the pullback from $C$ to $Q_d \times C$ of the line bundle ${\mathcal {O}}(m).$

The universal family $$\rho: 0 \rightarrow \e \rightarrow {\mathcal {O}}^N \rightarrow \f \rightarrow 0 \, \,
\, \text{on} \, \, \, Q_d\times C$$ gives rise to $$\eta_{\star} \rho (m): 0 \rightarrow
R^{0}\eta_{\star}\e (m) \rightarrow {\mathcal {O}}^N \otimes R^0 \eta_{\star}{\mathcal {O}}(m) \rightarrow
R^{0}\eta_{\star}\f (m) \rightarrow 0 \, \, \, \text{on} \, \, \, Q_d,$$ parametrizing vector
space quotients of a fixed vector space ${\bb C}^N \otimes {\bb C}^{\chi({\mathcal {O}}(m))}.$ This
association defines a closed embedding \begin{equation}\label{grothend} \iota: Q_d \rightarrow Grass(mr-d -r(g-1), Nm-
N(g-1)).\end{equation} For simplicity we denote this Grassmannian by $G$. On the level of closed points,
$\iota$ is realized as $$\left[0 \rightarrow E \rightarrow {\mathcal {O}}^N \rightarrow F \rightarrow
0\right] \mapsto \left[0 \rightarrow H^0(E(m)) \rightarrow {\bb C}^N \otimes {\bb
C}^{\chi({\mathcal {O}}(m))} \rightarrow H^0 (F(m)) \rightarrow 0\right].$$

Let $\k$ be the sheaf on $Q_d \times C$ defined by the exact sequence:
\begin{equation}
0 \rightarrow \k \rightarrow \eta^{\star}({\eta}_{\star} \e(m)) \otimes
 {{\mathcal {O}}}(-m) \rightarrow \e \rightarrow 0. \label{kernel} \end{equation} Applying the functor
$Hom_{\eta}(\cdot, \f)$ to \eqref{kernel} gives the sequence of sheaves on $Q_d$:
\begin{equation} 0 \rightarrow Hom_{\eta}(\e,\f) \rightarrow (\eta_{\star}\e(m))^{\vee} \otimes
\eta_{\star} \f(m) \rightarrow Hom_{\eta} (\k, \f) \rightarrow Ext^{1}_{\eta}(\e, \f)  
\rightarrow 0. \label{directimage} \end{equation} The second sheaf is computed by the projection
formula. It equals the pullback $\iota^{\star}\t G$ of the tangent bundle of $G$, so $\mathcal A_0:=\iota^{\star}\t G$.
Continuing the above exact sequence one more term, we obtain $Ext^{1}_{\eta}(\k, \f)=0$.
Therefore, the third sheaf $\mathcal A_1:=Hom_{\eta}(\k, \f)$ is locally free by cohomology and
base change.

To complete the derivation of $\eqref{tanobs}$ we observe that the argument goes through if we
twist by $\eta^{\star} \mathcal I$ {\it i.e.,} if we apply the functor $Hom_{\eta}(\cdot, \f\otimes
\eta^{\star}\mathcal I)$ to $\eqref{kernel}$. However, we do need to check that $$Hom_{\eta}(\k,
\f\otimes \eta^{\star} \mathcal I)= Hom_{\eta}(\k, \f)\otimes \mathcal I \;\;(=\mathcal
A_1\otimes \mathcal I).$$ The question is local in the base, so we may restrict to an affine open subscheme of $Q_d$. As in \cite {thomas}, the claim is a consequence of the spectral sequence below whose terms with $i$ or $j\geq 1$
vanish: $$\text{Tor}_{i}(Ext^{j}_{\eta} (\k, \f), \mathcal I)\Longrightarrow
\text {Ext}^{j-i}_{\eta}(\k,\f\otimes \eta^{\star}\mathcal I).$$ A similar argument holds for the other terms of $\eqref{tanobs}$.

Once the resolution of the tangent-obstruction complex $\eqref{tanobs}$ is obtained, the main
construction of \cite{litian} ensures the existence of a cone $\mathcal C$ inside the vector bundle
$\mathcal A_1$. The virtual fundamental class is standardly defined as the intersection of the class 
of this
cone with the zero section of $\mathcal A_1$: $$\left[Q_d\right]^{vir}=0^{!}_{\mathcal A_1}
\left[\mathcal C\right]\,\in\, H_{2e}(Q_d).$$ This concludes the proof of the theorem.

\begin {remark} Using the remarks following corollary $3.5$ in \cite {litian}, it should be possible to prove that $\mathcal C$ is the normal cone of the Grothendieck embedding $\eqref {grothend}$. There is another, admittedly clumsier method of obtaining an {\it equivalent} statement. Theorem $\ref{compatible}$ gives an alternate way of realizing the virtual fundamental class of the moduli space $Q_d$: it coincides with the virtual fundamental class of the zero locus of a section of a vector bundle over a local complete intersection; the latter is a Quot scheme for large degree $d$ having the expected dimension \cite {bertram1}. Using this new point of view, the construction can be recast in the more abstract language of cotangent complexes employed in \cite {behrendfantechi}. It is shown in \cite {pantev} that switching to the different formalism does not produce a new virtual class. One can now apply Siebert's work \cite {siebert} to write down the explicit formula: \begin {equation}\label{siebertclass} \left[Q_d\right]^{vir}=\left\{c(\mathcal A_1)\cap s(C_{Q/G})\right\}_{e}.\end{equation} Here, $s(C_{Q/G})$ is the Segre class of the normal cone for the Grothendieck embedding $\eqref{grothend}$, and $\left\{\cdot\right\}_{e}$ represents the component of degree equal to the virtual dimension $e$.\end{remark}

\begin {remark} In Gromov-Witten theory, one equips the scheme of degree $d$ morphisms to a Grassmannian $\mor$ with a virtual fundamental class \cite {behrendfantechi}. There is an open embedding \cite {bertram1} $$i: \mor\hookrightarrow Q_d,$$ alluded to in the introduction, which associates to a morphism $f:C\to G (r, N)$ the {\it locally free }quotient $$0\to f^{\star} S\to{\mathcal {O}}^N\to f^{\star} Q\to 0.$$ Here $$0\to S\to {\mathcal {O}}^N \to Q\to 0$$ is the tautological sequence on $G (r, N)$. It follows easily that $$i^{\star}[Q_d]^{vir}=[\text {Mor}_d]^{vir}.$$ \end {remark}

\begin{remark} Finally, we note that the construction can be done for the relative $Quot_{X/S}$ scheme
of quotients of any sheaf along the fibers of a smooth projective morphism $X\to S$ of relative
dimension $1$. Proposition $3.9$ in \cite {litian} immediately implies functoriality with respect to
cartesian diagrams. Nonetheless, for the purposes of this paper such level of generality will not be
needed.

\end {remark} 

\section{Equivariant localization}

In this section, we apply the virtual localization formula for $\quot$. We compute the fixed loci and their virtual normal bundles. 

\subsection {The virtual localization formula} Let us review the ingredients of the virtual
localization formula in \cite{graberpandharip}. To set the stage, let $X$ be a projective scheme over ${\bb
C}$ endowed with a $T={\bb C}^{\star}$ action, as well as with a perfect obstruction theory
compatible with this action. Let $X_i$ be the fixed loci of the action. The perfect obstruction
theory on $X$ induces perfect obstruction theories on each $X_i$: we consider only the fixed part
of the resolution $\eqref{tanobs}$. Similarly, the virtual normal bundle $\n_i^{vir}$ is the moving part of $\eqref{tanobs}$. If $j$ denotes the inclusion of the fixed loci into
$X$, the localization formula asserts that
\begin{equation} [X]^{vir} = j_{\star} \sum_{i} \frac{[X_i]^{vir}}{e_T(\n_i^{vir})}\, \, \,
\text{in} \, \, A_{T}(X) \otimes {\bb Q}\left[h, {h^{-1}}\right ].
\label{virtualloc} \end{equation} Here, we write $h$ for the generator of the equivariant
coefficient ring $H^{\star}_{T}(pt).$

The virtual localization formula derived in \cite {graberpandharip} uses the language of cotangent complexes and perfect obstruction theories of \cite {behrendfantechi}. Since the virtual localization formula was never written down in all generality using the Li-Tian formalism, we now briefly indicate why it holds in the particular situation of interest to us. While this is probably known to some experts, we include an outline of the argument for the convenience of the reader.

We have at our disposal an equivariant embedding
$$\iota: X \hookrightarrow M, \, \, M \,\text{smooth},$$ and a resolution of the tangent obstruction complex \eqref{tanobs} with ${\mathcal A}_0  = \iota^{\star} \t M.$ Let $M_i$ denote the fixed loci of the torus action on $M$, 
and denote by $\kappa_i: M_i \rightarrow M$ the inclusion. For each $i$, we have fiber square diagrams

\begin{center}
$\xymatrix{X_i \ar[r]^{j_i} \ar[d] & X \ar[d]
\\ M_i \ar[r]^{\kappa_i} &  M}$
\end{center}

Let $\m_i$ be the normal bundle of $M_i$ in $M$. The localization formula on 
$M$ gives
$$\sum_i (\kappa_i)_{\star} \frac{[M_i]}{e_{T}(\m_i)} = [M].$$ Taking refined intersections with $[X]^{vir}$ we have
$$\sum_i (j_i)_{\star} \frac{\kappa_i^{!} [X]^{vir}}{e_T(\m_i)} = [X]^{vir}.$$
Now the virtual normal bundle $\n_i^{vir}$ of $X_i$ in $X$ is $${\n}_i^{vir} = {\mathcal    
A}_0^m - {\mathcal A}_1^m = \m_i - {\mathcal A}_1^m.$$ Restrictions to $X_i$ are implied above. In addition, for a $T$-equivariant bundle on the fixed locus we use the superscripts ``$fix$" and ``$m$" to denote the fixed and moving subbundles. To prove \eqref{virtualloc} we only need to show
\begin{equation}
\label{jlocalize}
\kappa_i^{!} [X]^{vir} = e_T({\mathcal A}_1^m) \cap [X_i]^{vir}.\end{equation}
Now the perfect obstruction theory of $X_i$ as a fixed locus on $X$ is determined by the exact sequence
\begin{equation}\label{to1}0 \rightarrow \t^1 X_i \rightarrow {\mathcal A}_0^{fix} \rightarrow {\mathcal A}_1^{fix} \rightarrow \t^2 X_i\to 0 \text { and } {\mathcal A}_0^{fix} = \t M_i.\end{equation}
Note however that another valid perfect obstruction theory on $X_i$ is determined by 
\begin{equation}\label{to2}0 \rightarrow \t^1 X_i \rightarrow \t M_i \rightarrow {\mathcal A}_1 \rightarrow \widetilde{\t^2} X_i\to 0,\end{equation} where 
$\widetilde{\t^2} X_i = \t^2 X_i \oplus A_1^m.$ The virtual fundamental class $\widetilde {[X_i]}^{vir}$ of this second 
obstruction theory is related to $[X_i]^{vir},$ the virtual cycle of the first one, in an obvious manner:
\begin{equation}
\label{vir12}
\widetilde{[X_i]}^{vir} = e( {\mathcal A}_1^m) \cap [X_i]^{vir}.
\end{equation}
It is moreover this second perfect obstruction 
theory of $X_i$ which is compatible with the one on $X$ in the sense of Proposition 3.9 of \cite{litian}. The technical condition assumed there is satisfied: we have an 
exact sequence of complexes on $X_i$,
$$ 0 \rightarrow [0 \rightarrow \m_i] \rightarrow [\t M \rightarrow {\mathcal A}_1 \oplus \m_i] 
\rightarrow [\t M \rightarrow {\mathcal A}_1] \rightarrow 0,$$ whose cohomology yields, via \eqref{to1} and \eqref{to2}, the long exact sequence $$0\to \t^1 X_i \to \t^1 X \to \m_i|_{X_i} \to \widetilde \t^2 X_i \to \t^2 X\to 0.$$ Therefore, by \cite{litian} we have 
\begin{equation}
\label{functorial7}
\widetilde{[X_i]}^{vir} = \kappa_i^{!} [X]^{vir}. 
\end{equation}
Equations \eqref{vir12} and \eqref{functorial7} imply \eqref{jlocalize}, to which, it was argued, 
the virtual localization formula on $X$ is equivalent. 

Now, let $P(\alpha_1, \ldots, \alpha_n)$ be a polynomial in Chern classes of bundles on $X$, such that 
the degree of $P$ equals the virtual dimension of $X$, and such that $\alpha_1, \ldots, \alpha_n$ 
admit equivariant extensions $\tilde{\alpha}_1, \ldots, \tilde{\alpha}_n.$ Then 
\eqref{virtualloc} implies
\begin{equation}
\int_{[X]^{vir}} P(\alpha_1, \ldots, \alpha_n) = \sum_{i} \int_{[X_i]^{vir}} 
\frac{\iota^{\star} P(\tilde{\alpha}_1, \ldots, \tilde{\alpha}_n)}{e_T(\n_i^{vir})}, 
\label{virtualloc1}
\end{equation}
Our goal is to apply formula 
\eqref{virtualloc1} in the context of the Quot scheme for a suitable torus action.

\subsection{The torus action and its fixed loci}  Observe that $GL(N)$, viewed as the automorphism group of the trivial sheaf ${\mathcal {O}}^N$ on $C$, acts naturally on $\quot$: each quotient $\rho: {\mathcal {O}}^N \rightarrow F$ gets sent by $g\in 
GL(N)$ to $\rho \circ g: {\mathcal {O}}^N \rightarrow F$. We localize with respect to the 
diagonal subtorus $$T = \text{Diag}\, \{t^{-\lambda_1}, \ldots , t^{-\lambda_N} \}$$ 
of $GL(N),$ where the weights $\lambda_i, \, 1\leq i \leq N$ are distinct. We mention that the analysis of the fixed loci and of their normal bundles was also pursued in the case of the Quot scheme on ${\bb P}^{1}$ in \cite{strom}.

Since the $N$ copies of ${\mathcal {O}}$ in ${\mathcal {O}}^N$ are acted on with different weights, in order for a closed point
$$0 \rightarrow E \rightarrow {\mathcal {O}}^N \rightarrow F \rightarrow 0$$ to be fixed by $T$, it is
necessary and sufficient that $E$ split as a direct sum of line bundles, $$E=\oplus_{i=1}^{r}L_i,$$ each line bundle $L_i$ mapping
injectively to one of the $N$ available copies of ${\mathcal {O}}$. A fixed locus is thus labeled by a
choice of $r$ copies of ${\mathcal {O}}$ in ${\mathcal {O}}^N$ and the choice of an ordered $r$-partition of the degree
$d$, $d = d_1 + \cdots + d_r$. 

Having fixed these discrete data labeling a fixed locus $Z$, the point $$0 \rightarrow E
\rightarrow {\mathcal {O}}^N \rightarrow F \rightarrow 0$$ in $Z$ determines and is determined by $r$ exact
sequences \begin{equation}\label{exct}0 \rightarrow L_i \rightarrow {\mathcal {O}}_{k_i} \rightarrow T_i
\rightarrow 0, \, \, 1\leq i \leq r.\end {equation} Here the notation ${\mathcal {O}}_{k_i}$ is meant to
single out a copy of ${\mathcal {O}}$ from the available $N$, $\{k_1, \ldots , k_r \} \subset \{1, 2, \ldots,
N\},$ and the line bundle $L_i$ has degree $-d_i$. The fixed locus $Z$ corresponding to $\{k_1,
\ldots , k_r \} \subset \{1, 2, \ldots, N\}$ and to $\{d_1, \ldots , d_r \}$ is therefore a
product $$Z = Z_1 \times \cdots \times Z_r,$$ with $Z_i$ parametrizing exact sequences
$\eqref{exct}$, hence zero dimensional subschemes of $C$ of length $d_i$. Thus $Z_i$ is equal
to the $d_i^{\text {th}}$
symmetric product of $C$: 
\begin{equation}\label{fixedloc}Z_i \simeq \, {\text{Sym}}^{d_i} C, \,
\, \, \, \text{so} \, \, Z = {\text{Sym}}^{d_1} C \times \cdots \times {\text{Sym}}^{d_r} C.\end
{equation}

\subsection{The equivariant normal bundles}

We turn now to the problem of determining the weights of the $T$-action on the
virtual normal bundle to the fixed locus $Z$. For notational simplicity, we will assume that for quotients in $Z$, it is the first $r$ copies of ${\mathcal {O}}$ in ${\mathcal {O}}^N$ that are singled out, in other words that $k_{i} = i, \, 1 \leq i \leq r.$ The virtual normal
bundle to $Z$ is the {\em{moving}} part of the pullback to $Z$ of the virtual tangent bundle to $\quot,$ which can be read off from the resolution $\eqref{tanobs}$. According to
$\eqref{directimage}$, it equals
\begin{equation}
(\eta_{\star}\e (m))^{\vee} \otimes \eta_{\star} \f(m) - Ext^{0}_{\eta}(\k, \f).
\label{virtualtb}
\end{equation}  
We determine its restriction to $Z$. The universal sequence 
$$0 \rightarrow \e \rightarrow {\mathcal {O}}^N  \rightarrow \f \rightarrow 0$$ on $\quot 
\times C$ restricts to
\begin{equation}
0 \rightarrow {\mathcal {L}}_1 \oplus \cdots \oplus {\mathcal {L}}_r \rightarrow {\underbrace{{\mathcal {O}} \oplus 
\cdots \oplus {\mathcal {O}}}_{r}} \oplus
{{\mathcal {O}}}^{N-r}
\rightarrow \t_1 \oplus \cdots \oplus \t_r \oplus {{\mathcal {O}}}^{N-r} \rightarrow 0 
\label{univfixed}
\end{equation}
 on $Z 
\times C$. Here \begin{equation}\label{ustr}0 \rightarrow {\mathcal {L}}_i \rightarrow {\mathcal {O}} \rightarrow \t_i \rightarrow 
0\end{equation} on $Z_i \times C$ is
the universal structure associated with $Z_i \simeq {\text{Sym}}^{d_i} C$. For notational simplicity we denote by 
$0 \rightarrow {\mathcal {L}}_i \rightarrow {\mathcal {O}} \rightarrow \t_i \rightarrow 0$ not just the 
universal exact sequence on $Z_i \times C$, but also its pullback to $Z \times C,$
with $Z = Z_1 \times \cdots \times Z_r.$ 
The sheaf $\k$ on the other hand pulls back to 
\begin{equation}
\iota_{Z \times C}^{\star} \k = \oplus_{i=1}^{r}\k_i, 
\label{k}
\end{equation}
where $\k_i$ is defined by
\begin{equation}
0 \rightarrow \k_i \rightarrow \eta^{\star}(\eta_{\star} {\mathcal {L}}_i(m))\otimes
{\mathcal {O}}(-m) \rightarrow {\mathcal {L}}_i \rightarrow 0.
\label{ki}
\end{equation}
Equations \eqref{univfixed}, \eqref{k} and \eqref{ki} imply that the restriction of the 
virtual tangent bundle \eqref{virtualtb} to $Z$ is 
\begin{equation}
\label{messy}
\left(\bigoplus_{1\leq i, j \leq r} \left(\eta_{\star} {\mathcal {L}}_i (m)\right)^{\vee} \otimes \eta_{\star} 
\t_j (m) 
\mathop{\bigoplus_{1\leq i\leq r}}_{r+1\leq j\leq N}
\left(\eta_{\star} {\mathcal {L}}_i(m))^{\vee}
\otimes H^0 ({\mathcal {O}}(m)\right)_j \right)\ominus
\end{equation} 

$$\ominus\left(\bigoplus_{1 \leq i,j \leq r} 
\eta_{\star}\left(\k_i^{\vee} \otimes \t_j\right)
\mathop{\bigoplus_{1\leq i \leq r}}_{r+1 \leq j \leq N}
\left(\eta_{\star}\k_i^{\vee}\right)_j\right).$$ 
The subscript $j$ denotes tensor product with the $j^{\text {th}}$ copy of ${\mathcal {O}}$, and is used to keep track of the weights. From the ungainly expression {\eqref{messy}}, it is clear that the weights of the 
$T$-action on the restriction of the virtual tangent bundle {\eqref{virtualtb}} to $Z$ are as follows:

\begin{itemize} \item for $1\leq i, j\leq r$, the virtual sheaf \begin{equation}\label{nij1} \n_{ij} =
\left(\eta_{\star} {\mathcal {L}}_i (m)\right)^{\vee} \otimes \eta_{\star} \t_j(m) - \eta_{\star}\left(\k_i^{\vee}
\otimes \t_j\right)\end{equation} is acted on with weight $\lambda_i - \lambda_j$;\vskip.05in

\item for $1\leq i \leq r, \, r+1 \leq j \leq N$, the sheaf \begin{equation}\label{nij2}\n_{ij} =
(\eta_{\star} {\mathcal {L}}_i(m))^{\vee} \otimes H^0 ({\mathcal {O}}(m))_j - (\eta_{\star}\k^{\vee}_i)_j,\end{equation} is acted
on with weight $\lambda_i - \lambda_j$. \end{itemize}

In particular, the {\em{moving}} part of the virtual tangent bundle 
\eqref{virtualtb} 
on $Z$ is obtained by summing the nonzero-weight contributions in equations $\eqref{nij1}$ and $\eqref{nij2}$ above,
\begin{equation}\label{virtualnb}
\n^{vir} = \bigoplus_{\stackrel{1\leq i \leq r, 1 \leq j \leq N} {i \neq j}} \n_{ij}.
\end{equation}

\subsection{The equivariant Euler classes} We would like to further compute
the inverse of the equivariant Euler class appearing in equivariant localization $$\frac{1}{e_{T}(\n^{vir})}\in
H^{\star}(Z)\left[h, {h^{-1}}\right].$$

First, $\eqref{virtualnb}$ implies that
$$\frac{1}{e_{T}(\n^{vir})} = \mathop{\prod_{i,j}}_{i \neq j} \left ( \frac{1}{e_{T}(\n_{ij})} \right ).$$
We let $$c_t (\n_{ij}) = 1 + t c_1 (\n_{ij}) + t^2 c_2 (\n_{ij}) + \ldots$$ denote the total
Chern class of $\n_{ij}$, and we let $n_{ij}$ be the virtual rank of $\n_{ij}$. Since $\n_{ij}$
is acted on with weight $\lambda_i - \lambda_j$, we can furthermore write $$e_{T} (\n_{ij}) =
\left[(\lambda_i -\lambda_j)h\right]^{n_{ij}} c_{\tau_{ij}} (\n_{ij}),$$ where $$\tau_{ij} =
\frac{1}{(\lambda_i -\lambda_j)h}.$$
We conclude therefore that \begin{equation}\label{inverse}\frac{1}{e_{T}(\n^{vir})} =
\mathop{\prod_{1\leq i\leq r, 1\leq j\leq N}}_{i\neq j} \frac{1}{\left [(\lambda_i
-\lambda_j)h\right]^{n_{ij}}} s_{\tau_{ij}}(\n_{ij}),\end {equation} with $s_t (\n_{ij})$
denoting the total Segre class of $\n_{ij}.$

Note now that
\begin{itemize}
\item
for $1\leq i,j \leq r$, $i \neq j$, the sheaf $$\n_{ij} = (\eta_{\star} {\mathcal {L}}_i (m))^{\vee} \otimes 
\eta_{\star} \t_j(m) -
\eta_{\star}(\k_i^{\vee}
\otimes \t_j)$$ is $K$-theoretically equivalent to $\eta_{\star}({\mathcal {L}}_i^{\vee} \otimes \t_j).$ This 
can be 
seen from the defining sequence \eqref{ki} of $\k_i$, upon taking $Hom_{\eta}(\cdot \, , \t_j)$.
In order to compute $s_t(\n_{ij})$ we use the latter expression. Equation $\eqref{ustr}$ gives
\begin {equation}s_t(\n_{ij}) = s_t(\eta_{\star}({\mathcal {L}}_i^{\vee} \otimes \t_j))=s_t 
(\eta_{\star}({\mathcal {L}}_i^{\vee} - {\mathcal {L}}_i^{\vee} \otimes {\mathcal {L}}_j)) = s_t (\eta_{\star}({\mathcal {L}}_i^{\vee}))
\cdot c_t (\eta_{\star} ({\mathcal {L}}_i^{\vee} \otimes {\mathcal {L}}_j )).
\label{segre1}
\end {equation}
Moreover the virtual rank $n_{ij}$ of $\n_{ij}$ is 
\begin {equation}
n_{ij} = (d_i - \gbar) - (d_i - d_j - \gbar).
\label{rank1}
\end{equation}
In this discussion $\eta_{\star}=R^0\eta_{\star}-R^1\eta_{\star}$. \vskip.1in
\item
Similarly, for $1\leq i \leq r$ and $r+1 \leq j \leq N$, $$\n_{ij} = 
(\eta_{\star} {\mathcal {L}}_i(m))^{\vee} \otimes H^0 ({\mathcal {O}}(m))_j -
(\eta_{\star}\k^{\vee}_i)_j$$ is $K$-theoretically equivalent
to $\eta_{\star}({\mathcal {L}}_i^{\vee}).$ Therefore
\begin{equation}
s_t(\n_{ij}) = s_t (\eta_{\star}({\mathcal {L}}_i^{\vee})),
\label{segre2}
\end{equation}
and in this case
\begin{equation}
n_{ij} = d_i - \gbar.
\label{rank2}
\end{equation}
\end{itemize}

We express the Chern classes $c_t(\eta_{\star}({\mathcal {L}}_i^{\vee}))$ and $c_t (\eta_{\star}({\mathcal {L}}_i^{\vee} \otimes {\mathcal {L}}_j))$ in terms of classes on $Z_i$ and respectively $Z_i \times Z_j$, whose intersection theory is understood. Stated in general form the problem is as follows. Let $$0 \rightarrow {\mathcal {L}} \rightarrow {\mathcal {O}} \rightarrow \t \rightarrow 0$$ be the 
universal sequence on ${\text{Sym}}^{d} C \times C$. Write the K{\"{u}}nneth decomposition
of $c_1({\mathcal {L}}^{\vee})$ with respect to the chosen basis for $H^{\star} (C)$ as \begin {equation}\label{lvee}c_1({\mathcal {L}}^{\vee}) 
= x \otimes 1 + \sum_{1}^{2g} y^{j} \otimes \delta_j + d \otimes \omega.\end {equation}
It is well known \cite{arbarello} that $$(\sum_{j} y^{j} \otimes \delta_j)^2 = -2 \theta 
\otimes \omega,$$ where
$\theta$ denotes the pullback of the theta class under the map $${\text{Sym}}^{d} C 
\rightarrow \,\text{Pic}^{d}\,C.$$ Moreover, on ${\text{Sym}}^{d} C,$
\begin{equation}
\label{xtheta}
x^{d-l} \theta^{l} = \frac{g!}{(g-l)!}\, \, \text{for} \, \, l \leq g, \, \, \text{and} \, 
\, x^{d-l} \theta^{l} = 0 \text{ for} \, \,  l > g. 
\end{equation}

We compute $c_t(\eta_{\star}({\mathcal {L}}^{\vee}))$ in terms of $x$ and $\theta.$ This is
simple and amply explained in \cite{arbarello}. Indeed, by
the Grothendieck-Riemann-Roch theorem, one immediately finds $$\text{ch}\, (\eta_{\star}
({\mathcal {L}}^{\vee})) = e^{x} (d -\gbar - \theta).$$ Quite generally, the total Chern class of a rank $n$ vector
bundle $H$ with Chern character $ n -\theta$ is $e^{-t\theta}.$ The total Chern class of $H
\otimes L$, where $L$ is a line bundle with first Chern class $x$, is $$c_t (H \otimes L) =
(1+tx)^n c_{\tau} (H)=(1+tx)^{n}e^{-\tau \theta}\,, \, \, \text{with} \, \, \tau = \frac{t}{1+tx}.$$ This observation gives
the Chern and Segre polynomials:
\begin{equation}
s_t(\eta_{\star}({\mathcal {L}}^{\vee})) = \frac{1}{c_t(\eta_{\star}({\mathcal {L}}^{\vee}))}=(1+tx)^{-d+\gbar}
\exp\left({\frac{t}{1+tx} \theta}\right). \label{equivsegre1} \end{equation}

By the same token, $$\text{ch}\, (\eta_{\star} ({\mathcal {L}}_i^{\vee}\otimes {\mathcal {L}}_j)) = e^{x_i - x_j} (d_i - 
d_j - \gbar -
(\theta_i + \theta_j + \sigma_{ij})),$$ with \begin {equation}\sigma_{ij} = -
\sum_{k=1}^{g}
\left(y_i^{k} y_j^{k+g} + y_j^{k} y_i^{k+g}\right).\end {equation}
Therefore, \begin{equation} 
c_t(\eta_{\star}({\mathcal {L}}_i^{\vee}
\otimes {\mathcal {L}}_j)) = (1+t(x_i - x_j))^{d_i- d_j -\gbar} \exp \left({-\frac{t}{1+t(x_i - x_j)} (\theta_i +
\theta_j + \sigma_{ij})}\right). \label{equivchern2} \end{equation} 

Together, \eqref{inverse}, \eqref{segre1}, \eqref{rank1}, \eqref{segre2}, \eqref{rank2},
\eqref{equivsegre1}, and \eqref{equivchern2} give the indecorous result
\begin{eqnarray*}
\frac{1}{e_{T}(\n^{vir})} &=& {\mathop{\prod_{1\leq i, j \leq
r}}_{i\neq
j}}
\tau_{ij}^{d_j - d_i + \gbar}
\left ( 1+ \tau_{ij} (x_i - x_j )\right )^{d_i - d_j - \gbar} 
\exp\left({-\frac{\tau_{ij}}{1+\tau_{ij}(x_i - x_j)}
(\theta_i + \theta_j + \sigma_{ij})}\right) \\ & \cdot &
{\mathop{\mathop{\prod_{1\leq i \leq r}}_{1\leq j \leq N}}_{i\neq j}} 
\tau_{ij}^{d_i 
- \gbar} (1+\tau_{ij} x_i)^{-d_i + \gbar} 
\exp\left({\frac{\tau_{ij}}{1+\tau_{ij}x_i} \theta_i}\right),
\end{eqnarray*}
Since $\tau_{ij} = -\tau_{ji},$
the first product simplifies greatly, as it ranges over all pairs $(i,j), \,1\leq i 
\neq j\leq r$. Thus, we can rewrite 
\begin{equation}\label{prel}
\frac{1}{e_{T}(\n^{vir})}=(-1)^{\gbar \binom{r}{2}+d(r-1)}{\mathop{\prod_{1\leq i<j \leq r}}}
\tau_{ij}^{2\gbar}
\left (1+ \tau_{ij} (x_i - x_j)\right )^{-2\gbar}\cdot\end{equation}
$$\cdot \prod_{1\leq i \leq r}\mathop{\prod_{1\leq j \leq N}}_{j\neq i} 
\tau_{ij}^{d_i - \gbar} (1+\tau_{ij} x_i)^{-d_i + \gbar} 
\exp\left({\frac{\tau_{ij}}{1+\tau_{ij}x_i} \theta_i}\right)
$$

The result of any localization computation does not depend on the specific weights that one uses. For the rest of the paper, we choose the $\lambda_i$s to 
be the $N^{\text {th}}$ roots of unity $$\lambda_i=\exp\left(\frac{2\pi i \sqrt{-1}}{N}\right), \, \, 1\leq i \leq N.$$ With this choice, the second 
product 
in $\eqref{prel}$ simplifies, since:
$$ {\mathop{\prod_{j}}_{j\neq i}} \left ( \frac{1}{\tau_{ij}} + x_i \right ) = 
{\mathop{\prod_{j}}_{j\neq i}} \left ( (\lambda_i h + x_i) - \lambda_j h \right ) = 
\frac{(\lambda_i h + x_i)^N - h^N}{x_i}, \, \, \, \, \text{and}$$
$${\mathop{\sum_{j}}_{j\neq i}} \frac{1}{{\tau_{ij}}^{-1} + x_i} = {\mathop{\sum_{j}}_{j\neq 
i}} \frac{1}{(\lambda_i h +x_i) - \lambda_j h} = \frac{N (\lambda_i h + x_i)^{N-1}}{(\lambda_i h + 
x_i)^N - h^N} - \frac{1}{x_i}.$$
We obtain the following expression for the Euler class,

\begin{equation}\label{euler2}
\frac{1}{e_{T}(\n^{vir})} = {\prod_{1\leq i \leq r}}
\left ( \frac{(\lambda_i h + x_i)^N - h^N}{x_i} \right )^{-d_i + \gbar}
\exp\left(\theta_i\cdot{\left ( \frac{N (\lambda_i h + x_i)^{N-1}}{(\lambda_i h +
x_i)^N - h^N} - \frac{1}{x_i}\right )}\right)\end {equation}

$$\cdot{\mathop{\prod_{1\leq i < j \leq r}}}
\tau_{ij}^{2\gbar}
\left (1+ \tau_{ij} (x_i - x_j)\right )^{-2\gbar}\cdot (-1)^{\gbar \binom{r}{2}+d(r-1)}.
$$

Finally, we remark that the argument leading to equations $\eqref{segre1}$ and $\eqref{rank1}$ shows that the induced 
virtual tangent bundle to the fixed locus $Z$ is $K$-theoretically equivalent to the actual tangent bundle of $Z$.
Hence the induced virtual fundamental class on the fixed locus $Z=\text {Sym}^{d_1} C\times \ldots \text {Sym}^{d_r} C$ coincides 
with the usual fundamental class.

\subsection {Compatibility of the virtual classes} In this subsection, we will use the localization computations above to
give a proof of Theorem $\ref{compatible}$. 

We let $$i_p: Q_d \to Q_{d+r}$$ be the equivariant embedding of Quot schemes which on closed points is given by
$$\left[E\hookrightarrow \mathcal O^N \right] \mapsto \left[\tilde E=E(-p)\hookrightarrow \mathcal O^N(-p)\hookrightarrow
\mathcal O^N\right].$$ We denote by $\mathcal {\tilde E}$ the universal subsheaf on $\tilde Q=Q_{d+r}$; clearly
$i_p^{\star} \mathcal {\tilde E}=\mathcal E(-p)$. The dual $$\left(\mathcal O^N_{p}\right)^{\vee} \to \mathcal {\tilde
E}_{p}^{\vee}$$ of the natural inclusion vanishes on $Q_d$. Just as in the non-virtual case, we seek to show that
\begin{equation}\label{vir0}i_{\star}\left[Q_d\right]^{vir}=c_{\text {top}}(Hom (\tilde {\mathcal E}_p, \mathcal
O^N))\cap \left[Q_{d+r}\right]^{vir}.\end {equation} We will establish $\eqref{vir0}$ equivariantly; we regard $\mathcal
O^N$ as an equivariant sheaf with the usual weights. We compute both sides by the localization formula $\eqref{virtualloc}$, and match the fixed loci contributions. 

Each fixed locus on $Q_d$, $$Z=\text {Sym}^{d_1} C\times \ldots \times \text {Sym}^{d_r}C,$$ corresponds to a fixed locus
on $\tilde Q,$ $$\tilde Z= \text {Sym}^{d_1+1} C \times \ldots \times \text {Sym}^{d_r+1} C.$$ There is an inclusion 
$j_p:
Z\to \tilde Z$ given on closed points by 
$$[L_i \hookrightarrow {\mathcal O}] \mapsto [\tilde L_i = L_i(-p) \hookrightarrow {\mathcal O}(-p) \hookrightarrow 
{\mathcal O}], \, \, \, 1\leq i \leq r.$$ Let $\tilde {\mathcal L}_i$ denote the pullback to $\tilde Z\times C$ of the universal bundle on the $i^{\text{th}}$ factor, and set $\tilde {x}_i=c_1(\tilde {\mathcal L}_i^{\vee}|_p)$. It is clear that there is a section of the bundle 
$\bigoplus_{i=1}^{r} 
{{{\tilde{{\mathcal {L}}}^{\vee}}_i|_{p}}}$ on $\tilde Z$ vanishing precisely on $Z$, such that \begin 
{equation}\label{vir1} ({j_p})_{\star}[Z]=\prod_{i=1}^{r}
\tilde x_i \cap [\tilde Z].\end {equation} Next, the following equivariant equality holds on $\tilde Z$: \begin{equation}\label{vir3}c_{\text {top}}(Hom(\mathcal {\tilde
E}_p, \mathcal O^N))_{\vert \tilde Z}=\prod_{i=1}^{r} \prod_{j=1}^{N} (\lambda_i h + \tilde x_i - \lambda_j
h)=\prod_{i=1}^{r} \tilde x_i \cdot \prod_{i=1}^{r}\left(\frac{(\lambda_i h + \tilde x_i)^N - h^N}{\tilde x_i} 
\right).\end{equation}
Moreover, using $\eqref{euler2}$ we obtain \begin{equation}\label{vir2}\frac{1}{e_{T}(\mathcal
N^{vir}_{Z/Q_d})}=j_p^{\star} \prod_{i=1}^{r}\left(\frac{(\lambda_i h + \tilde x_i)^N - h^N}{\tilde x_i} \right)\cdot 
j_{p}^{\star}
\frac{1}{e_{T}(\mathcal N^{vir}_{{\tilde Z}/{\tilde Q}})}.\end {equation} The localization contributions from $Z$ and $\tilde
Z$ to the two sides of $\eqref{vir0}$ are now matched using equations $\eqref{vir1}$, $\eqref{vir3}$, $\eqref{vir2}.$

Finally, there are other fixed loci on $\tilde Q$ that we need to consider. However, their contribution to the terms in
$\eqref{virtualloc1}$ vanishes. Indeed, for the remaining fixed loci it must be that one of the degrees $\tilde d_i$ is 
$0$.
The contribution of $c_{\text {top}} (Hom(\mathcal {\tilde E}_p, \mathcal O^N))$ on $\text {Sym}^{0}C$ must vanish
because of the term $\tilde x_i$ in expression $\eqref{vir3}$. This completes the proof of Theorem
$\ref{compatible}$.

\section{The Vafa-Intriligator formula} In this section, we derive the formula of Vafa-Intriligator for the intersections
of $a$-classes. The essential part of the argument is the {\it summation} of the fixed loci contributions which were
computed in the previous section.
 
To begin, let $P (a_1, \ldots, a_r)$ be a polynomial in $a$-classes of degree equal to the expected
dimension $e$ of $Q_d$, $ e = Nd -r (N-r)\gbar.$ We would like to compute $$P (a_1, \ldots, a_r)\cap [Q_d]^{vir}.$$
Write $z_1, \ldots, z_r$ for the Chern roots of $\e_p^{\vee},$ so $$P (a_1, \ldots, a_r) =
R (z_1, \ldots, z_r)$$ for a polynomial $R$. Explicitly, if $\sigma_1, \ldots, \sigma_r$ denote the elementary symmetric functions, we have $$R(z_1, \ldots, z_r)=P(\sigma_1(z_1, \ldots, z_r), \ldots, \sigma_r(z_1, \ldots, z_r)).$$

We first examine the pullbacks of the equivariant $a_i$s to the fixed locus $Z$ determined 
by the choice of an ordered degree splitting $d= d_1 + \cdots + d_r$ and of $r$ roots of unity
$\lambda_1, \ldots, \lambda_r$ out 
of the available $N$; each root of unity represents the weight on one of the $N$ copies of ${\mathcal {O}}$. It was already noted that $$\iota_{Z \times C}^{\star} \e = {\mathcal {L}}_1 \oplus \ldots \oplus {\mathcal {L}}_r.$$ 
Moreover,
using the notation of equation \eqref{lvee}, the equivariant Chern class $c_1^{T} ({\mathcal {L}}_l^{\vee})$ is 
\begin{equation}
\label{equivli}
c_1^{T} ({\mathcal {L}}_l^{\vee}) = x_l \otimes 1 + \sum_{j=1}^{2g} y_l^j \otimes \delta_j + d_l
\otimes \omega + \lambda_l h, \, \, \, 1\leq l \leq r.
\end{equation}
Therefore, we have
$$\iota_Z^{\star} a_i = \sigma_i \left ( (\lambda_1 h + x_1), \ldots, (\lambda_r h + x_r) \right ).$$ Moreover, \begin {equation} i_{Z}^{\star} P(a_1, \ldots, a_r)=R \left ( (\lambda_1 h + x_1), \ldots, (\lambda_r h + x_r) \right).\end{equation}

Let us write ${Z_{\{d_i\}}} = \text{Sym}^{d_1} C\times \cdots \times \text{Sym}^{d_r}C$. It is convenient to set \begin{equation}u=(-1)^{\gbar\binom{r}{2}+d(r-1)}.\end{equation} To prove Theorem $\ref{vafaint}$ using the virtual 
localization formula, it is therefore enough to show:
\begin{lemma}\label{sumlemma}
Fix any distinct $N$-roots of unity $\lambda_1, \ldots, \lambda_r,$ corresponding to a choice of rank 
$r$ trivial subbundle of ${\mathcal {O}}^N$. Let $R (z_1, \ldots, z_r) = z_1^{\alpha_1}\cdots z_r^{\alpha_r}$ be a 
monomial of degree $e$. Then  
\begin{equation}
\sum_{d_1 + \cdots d_r = d} \int_{Z_{\{d_i\}}} \frac{R(\lambda_1 h + x_1, 
\ldots, 
\lambda_r h + x_r)}{e_{T} (\n^{vir}_{Z_{\{d_i\}}})} = u\cdot R(\lambda_1, \ldots, 
\lambda_r) J^{g-1} (\lambda_1, 
\ldots, \lambda_r)
\end{equation}
where $$J(\lambda_1, \ldots, \lambda_r) = N^r \cdot {\lambda_1^{-1} \cdot\ldots\cdot \lambda_r^{-1}}\cdot {\prod_{i<j} 
(\lambda_i - \lambda_j)^{-2}}.$$ 
\end{lemma}

{\it{Proof.}} By \eqref{euler2}, the integral $$\int_{Z_{\{d_i\}}} \frac{R(\lambda_1 h + x_1,
\ldots,
\lambda_r h + x_r)}{e_{T} (\n^{vir}_{Z_{\{d_i\}}})}$$ involves only $x$ and $\theta$ classes, whose 
intersection theory is given by \eqref{xtheta}. The $\theta$'s only appear in the exponentials 
\begin{equation}
\label{exponential}
\exp\left(\theta_i\cdot {\left ( \frac{N (\lambda_i h + x_i)^{N-1}}{(\lambda_i h +
x_i)^N - h^N} - \frac{1}{x_i}\right )}\right) = \sum_{l=0}^{g} \frac{{\theta_i}^l}{l!} 
\left(\frac{N (\lambda_i h + x_i)^{N-1}}{(\lambda_i h +
x_i)^N - h^N} - \frac{1}{x_i}\right )^l.
\end{equation}
For the purpose of intersecting
$x$ and $\theta$ classes on ${Z_{\{d_i\}}} = \text{Sym}^{d_1} C\times \cdots \times 
\text{Sym}^{d_r}C,$ we can replace with impunity 
$\theta_i^{l}$ by $\frac{g!}{(g-l)!} x_i^l,$ and the 
exponential \eqref{exponential} by 
\begin{eqnarray*}
\label{exp1}
\sum_{l=0}^{g} \frac{1}{l!} \frac{g!}{(g-l)!} x_i^l
{\left ( \frac{N (\lambda_i h + x_i)^{N-1}}{(\lambda_i h +
x_i)^N - h^N} - \frac{1}{x_i}\right )}^l &= &\left ( 1 + x_i 
\left ( \frac{N (\lambda_i h + x_i)^{N-1}}{(\lambda_i h +
x_i)^N - h^N} - \frac{1}{x_i}\right ) \right )^{g}\\ & = &N^g {x_i}^g\frac{ (\lambda_i h + 
x_i)^{(N-1)g}}{{\left ( (\lambda_i h + x_i)^N -h^N \right )}^{g}}. 
\end{eqnarray*}
Thus by \eqref{euler2},
\begin{equation}
\label{conc1}
\int_{Z_{\{d_i\}}} \frac{R(\lambda_1 h + x_1,
\ldots, \lambda_r h + x_r)}{e_{T} (\n^{vir}_{Z_{\{d_i\}}})} = u N^{rg}\cdot \,I_{\{d_i\}}\,\cdot h^{-\gbar r(r-1)} {\mathop{\prod_{1\leq i<j \leq r}}} (\lambda_i- \lambda_j)^{-2\gbar} 
\end {equation}
where we define
\begin{equation}\label{idi}I_{\{d_i\}}=\int_{Z_{\{d_i\}}} \prod_{1\leq i \leq r} 
(\lambda_i h + 
x_i)^{\alpha_i + (N-1)g} \frac{x_i^{d_i +1}}{\left ( (\lambda_i h + x_i)^N -h^N \right )^{d_i 
+1}}
\cdot {\mathop{\prod_{1\leq i< j \leq r}}}
\left (1+ \tau_{ij} (x_i - x_j)\right )^{-2\gbar} .\end {equation}
Since $$\int_{\text{Sym}^{d_i} C} x_i^{d_i} = 1,$$ $I_{\{d_i \}}$ is the same as the coefficient of 
$x_1^{d_1} \cdots x_r^{d_r}$ in its defining integrand.
Factoring out the $\lambda_i h$'s and letting formally
$$\bar{x}_i = \frac{x_i}{\lambda_i h},$$ we have equivalently
\begin{equation}
\label{residue3}
h^{-\gbar r (r-1)} I_{\{d_i \}} = R(\lambda_1, \ldots, \lambda_r) \cdot \prod_{1 \leq i \leq r} \lambda_i^{-\gbar} 
\cdot\bar{I}_{\{d_i \}},
\end{equation}
where we put 
\begin{equation}
\bar{I}_{\{d_i \}} = \text{Res }_{\bar{x}_i = 0}\left ( \prod_{1\leq i \leq r} 
\frac{(1+\bar{x}_i)^{\alpha_i + 
(N-1)g}}
{\left ( (1 + \bar{x}_i)^N -1 \right )^{d_i +1}} \prod_{1\leq i <j \leq r} \left ( 1 + \frac{\lambda_i 
\bar{x}_i - \lambda_j \bar{x}_j}{\lambda_i - \lambda_j} \right )^{-2\gbar} \right ).
\label{residue1}
\end{equation}

The following lemma would lead to an explicit expression for this residue
were it not for the presence of the {\it mixed} term $\prod_{1\leq i <j \leq r} \left 
( 1 + 
\frac{\lambda_i
\bar{x}_i - \lambda_j \bar{x}_j}{\lambda_i - \lambda_j} \right )^{-2\gbar}$ which renders the 
combinatorics slightly unpleasant.

\begin {lemma}\label{lemmaresidue}
The residue of $$\frac{(1+x)^{N-1+l}}{\left((1+x)^{N}-1)\right)^{d+1}}\cdot x^m$$  at $x=0$ is computed by the binomial sum $$\frac{1}{N}\sum_{p=0}^{m}(-1)^{m-p} \binom{m}{p}\binom{\frac{l+p}{N}}{d}.$$
\end {lemma}

{\it Proof.} The case $m=0$ is equivalent to the equation \begin {equation}\label{residue4}\text
{Res}_{x=0}\;\frac{(1+x)^{N-1+l}}{((1+x)^{N}-1)^{d+1}}=\frac{1}{N} \binom{\frac{l}{N}}{d}.\end {equation} The statement
for arbitrary $m$ follows from this case and from the identity $$x^m=\sum_{p=0}^{m} (-1)^{m-p} \binom{m}{p} (1+x)^p.$$ 

To prove $\eqref{residue4}$ we note that both sides are polynomials in $l$, hence we only need to establish the equality for infinitely many values of $l$, such as all multiples of $N$. Writing $$f_{l,
d} (x) =\frac{(x+1)^{N-1+Nl}}{((x+1)^{N}-1)^{d+1}}$$ we will show \begin {equation}\label{residue2}\text {Res}_{x=0}
\;f_{l,d}(x)=\frac{1}{N}\binom{l}{d}.\end{equation} We observe that $$f_{l+1, d}=f_{l, d-1}+f_{l,d}.$$ Our equation
$\eqref{residue2}$ follows inductively. We need to check the base cases $d=0$ and $l=0$. Only the second check requires
an explanation; in this case, $$\text {Res}_{x=0} \; \frac{(x+1)^{N-1}}{\left((x+1)^{N}-1\right)^{d+1}}= - 
\frac{1}{Nd} \text { Res}_{x=0}
\;\partial_{x} \left((x+1)^N-1\right)^{-d}= 0,$$ as desired. The last equality is a general fact about residues of
derivatives of meromorphic functions. Our proof is complete. \\

In order to apply the lemma to \eqref{residue1}, we write
\begin{eqnarray*}
&\prod_{1\leq i <j \leq r} &\left
( 1 + \frac{\lambda_i \bar{x}_i - \lambda_j \bar{x}_j}{\lambda_i - \lambda_j} \right )^{-2\gbar}
 = \sum_{k_{ij}\geq 0}\, \prod_{i<j} \binom{-2\gbar}{k_{ij}}\left ( \frac{\lambda_i 
\bar{x}_i - \lambda_j 
\bar{x}_j}{\lambda_i - \lambda_j} \right )^{k_{ij}}=\\
& = & \sum_{k_{ij}\geq 0}\, \sum_{m_{ij} = 
0}^{k_{ij}} \prod_{i<j} (-1)^{k_{ij} - m_{ij}} \binom{-2\gbar}{k_{ij}}\binom{k_{ij}}{m_{ij}} 
\frac{(\lambda_i
\bar{x}_i)^{m_{ij}} (\lambda_j
\bar{x}_j)^{k_{ij} - m_{ij}}}{(\lambda_i - \lambda_j)^{k_{ij}}}.
\end{eqnarray*} 
In this summation we allow all possible level values $$k = \sum_{i < j} k_{ij}, \, \, \,k 
\leq d.$$
For $j > i$ we define $m_{ji}  = k_{ij} - m_{ij}.$ Then $\bar x_i$ appears in the above sum
with the power $m_i =_{\text{def}} \sum_{j\neq i} m_{ij},$ and we also have $\sum_{i} m_i = k.$

We use Lemma \ref{lemmaresidue} to determine the coefficient of \begin {equation}\label{coef}\prod_{i<j} (-1)^{k_{ij} - m_{ij}}
\binom{-2\gbar}{k_{ij}}\binom{k_{ij}}{m_{ij}} \frac{\lambda_{i}^{m_{ij}} \lambda_j^{k_{ij} -
m_{ij}}}{(\lambda_i-\lambda_j)^{k_{ij}}}\end {equation} in the {\it{sum over degrees}} $$\sum_{d_1 + \cdots + d_r = d}
\bar{I}_{\{d_i \}}.$$ Set $l_i = \alpha_i + (N-1)\gbar,$ so $$\sum_{i} l_i= Nd+ \gbar r(r-1)=_{\text {def}}
Nd+s.$$ The coefficient in question equals
$$\frac{1}{N^r}\sum_{d_i} \prod_{i=1}^{r} \sum_{p_i = 0}^{m_i} (-1)^{m_i-p_i} 
\binom{m_i}{p_i} \binom{\frac{l_i + p_i}{N}}{d_i} = \frac{1}{N^r}\prod_{i=1}^{r} (-1)^{m_i} \sum_{p_i = 
0}^{m_i}
(-1)^{p_i} \binom{m_i}{p_i} \binom{\frac{\sum_{i} l_i + \sum_{i} p_i}{N}}{d}$$
$$=\frac{1}{N^r}(-1)^k \sum_{n=0}^{k} (-1)^n \binom{k}{n} 
\binom{d+\frac{s+n}{N}}{d}.$$
Here, we first summed over the degrees $\sum_{i} d_i=d$, then we summed over $\sum_{i} p_i = n$, recalling that $\sum {m_i} = k$. 

Finally, after summing over $m_{ij}\;$, we conclude that $$\sum_{d_1 + \cdots + d_r = d} \bar{I}_{\{d_i \}}=\frac{1}{N^r}\sum_{k=0}^{d} (-1)^k \sum_{n=0}^{k} (-1)^n \binom{k}{n}   
\binom{d+\frac{s + n}{N}}{d} \cdot \sum_{\sum k_{ij} =k} \prod_{i <j} 
\binom{-2\gbar}{k_{ij}}.$$
This simplifies to 
\begin{equation}
\label{almostfinal}
\sum_{d_1 + \cdots + d_r = d}\bar{I}_{\{d_i \}} = \frac{1}{N^r}\sum_{k=0}^{d} \sum_{n=0}^{k} 
(-1)^{n+k} 
\binom{k}{n}
\binom{d+\frac{s + n}{N}}{d} \binom{-\gbar r(r-1)}{k}=\frac{1}{N^r}.
\end{equation}
The last equality is the content of the next lemma. Together, equations $\eqref{conc1}$, $\eqref{residue3}$, 
and $\eqref{almostfinal}$ conclude the proof of Lemma $\ref{sumlemma}$.  

\begin {lemma} Let $d, N, s\geq 0$ be integers. Then $$\sum_{n=0}^{d}\sum_{k=0}^{d}(-1)^{k-n}\binom{-s}{k}\binom{k}{n}
\binom{d+\frac{s+n}{N}}{d}=1.$$
\end {lemma}

{\it Proof.} Let us define the integers $$c_n=\sum_{k=0}^{d}(-1)^{k-n}\binom{-s}{k}\binom{k}{n},\;\;\;\; 0\leq n\leq d.$$
We will show $$\sum_{n=0}^{d} c_n \binom {d+\frac{s+n}{N}}{d}=1.$$

We will work in the algebra $A=\mathbb C[x]/(x^{d+1})$. It is clear that for any nilpotent element $u\in A$, we can
define the nilpotent element $$\log(1+u)=-\sum_{i\geq 1} \frac{u^{i}}{i}.$$ Exponentiating nilpotents causes no
convergence problems; this allows us to define all powers $$(1+u)^{\alpha}=\exp (\alpha \log(1+u)).$$ The usual binomial
formulas are true for formal reasons.

Let us define $$y=1-(1+x)^{-N}.$$ The expansion of $y$ as a polynomial in $x$ has no constant term so $y^{d+1}=0$. It is
clear that we can solve for $$x=-1+(1-y)^{-\frac{1}{N}}$$ as a polynomial in $y$. This implies that the elements $1, y,
\ldots, y^{d}$ span the vector space $A$. Hence, they must form a basis for $A$.

We note that \begin {equation}\label{sumcp}\sum_{n=0}^{d} c_n(1+x)^{n}=\sum_{k=0}^{d} \binom{-s}{k} \sum_{n=0}^{k}
(-1)^{k-n}\binom{k}{n}(1+x)^{n}=\sum_{k=0}^{d}\binom{-s}{k} x^k = \frac{1}{(1+x)^{s}}.\end {equation} Using $\eqref{sumcp}$ and the standard binomial identities, we observe that $$\sum_{j=0}^{d}\left(\sum_{n=0}^{d}
c_{n}\binom {j+\frac{s+n}{N}}{j}\right) y^j = \sum_{n=0}^{d}\sum_{j=0}^{d} c_n \binom{-\frac{s+n}{N}-1}{j} (-y)^{j}=
\sum_{n=0}^{d} c_n (1-y)^{-\frac{s+n}{N}-1}$$ $$= \frac{1}{1-y} \sum_{n=0}^{d} c_n (1+x)^{s+n} =
\frac{1}{1-y}=\sum_{j=0}^{d} y^j.$$ The conclusion of the lemma now follows by considering the coefficient of $y^d$.

\section {Other intersection numbers} 

The intersection theory of $a$-classes is not sufficient for all applications. For example, if one wishes to compute the
degree of the Grothendieck embedding $\eqref{grothend}$ of $\quot$ into the Grassmannian, more general intersections need
to be calculated.

It is clear that in principle our method extends to compute all virtual intersection numbers on $\quot$. We use 
the virtual localization formula $\eqref{virtualloc1}$, and the computation of the virtual Euler classes in
$\eqref{euler2}$. In addition, we make use of the intersection theory of symmetric products. Note that the rules of
\eqref{xtheta} do not a priori determine all of the intersection numbers on $\text {Sym}^d C$ that we need. The
intersections of $x$ and $\theta$ with a product of $y^j$s are specified by the following observations, which can be
proved following the arguments of \cite {thaddeus}. They are sufficient for the proof of Proposition $\ref {allint}$.  
\begin{enumerate} \item [(i)] In a nonzero top intersection product on ${\text{Sym}}^d C,$ $y^j$ appears if and only if
$y^{j+g}$ appears as well. Since these classes are odd, they should appear only with exponent $1$.  \item [(ii)] The
product $$y^{j_1} y^{j_1 + g} \cdots y^{j_n} y^{j_n + g}\cdot P(x, \theta) \cap [{{\text{Sym}}^{d} C}]$$ is independent
of the specific distinct $j_1, \ldots, j_n \in \{1, \ldots , g\}.$ Thus, since $\theta = \sum_{j=1}^{g}y^j y^{j+g},$ this
intersection product equals $$\frac{(g-n)!}{g!} \theta^n \cdot P(x, \theta) \cap [{{\text{Sym}}^d C}].$$ \end{enumerate}

It remains to determine the equivariant restrictions of the $b$ and $f$ classes to the fixed locus $Z_{\{d_i\}}$. This is achieved in the formulas below. We let \begin{itemize}
\item $\sigma_i (x_1, \ldots, x_r)$ denote the $i^{\text{th}}$
symmetric function in $x_1, \ldots, x_r,$ 
\item $\sigma_{i;k} (x_1,
\ldots, x_r)$ be the $i^{\text {th}}$ symmetric function in the $r-1$ variables
$\{x_1 \ldots x_r\} \setminus \{x_k \}$, 
\item $\sigma_{i; k,l}$ be the
$i^{\text {th}}$ symmetric function in the $r-2$ variables $\{x_1, \ldots, x_r\} \setminus
\{x_k, x_l\}.$\end {itemize}

Recall that $\iota_{Z \times C}^{\star}\e = {\mathcal {L}}_1 \oplus \cdots \oplus
{\mathcal {L}}_r.$ Taking account of equation \eqref{equivli}, and setting $$\tilde x_i= x_i+\lambda_i h,$$
we get the following expressions of the 
equivariant restrictions of $a$, $b$ and $f$ classes:

\begin{itemize}
\item
$a_i = \sigma_i (\tilde x_1, \ldots, \tilde x_r)$
\item
$b_i^j = \sum_{q=1}^{r} y_q^j \;\sigma_{i-1;q} (\tilde x_1, 
\ldots , \tilde x_r)$
\item
$f_i = - \sum_{j=1}^{g} \sum_{a} \sum_{b\neq a} y_a^j y_b^{j+g}\;
\sigma_{i-2;a, b}(\tilde x_1, \ldots, \tilde x_r) + \sum_{q=1}^{r} d_q \sigma_{i-1;q}(\tilde x_1, \ldots, \tilde x_r)$.
\end{itemize}

For example, 
$$f_2=- \sum_{j=1}^{g}\mathop{\sum_{1\leq a, b\leq r}}_{a\neq b} y_a^{j} y_{b}^{j+g} + \sum_{q=1}^{r}
d_q \sum_{i\neq q} (x_i+\lambda_i h).$$

Intersections involving only $a$ and $b$ classes are explicitly computable with the methods of this work. For instance,
Proposition \ref {bees} gives all the required intersections in rank $2$. \\

{\it Proof of Proposition $\ref {bees}$.} Fixing $s\leq d,$ $s \leq g,$ and $1\leq j_1 
<\ldots
<j_s\leq g$, we show $$\left(b_1^{j_1}b_1^{j_1+g} \cdot \ldots\cdot b_1^{j_s} b_1^{j_s+g}\right)\cdot P(a)\cap \left[Q_d\right]^{vir} =\frac{u} {N^{s}} \sum_{\lambda_1, \ldots, \lambda_r} (\lambda_1+\ldots+\lambda_r)^s
\cdot \left(R J^{\gbar}\right)(\lambda_1, \ldots, \lambda_r).$$ Using the previous observations $(i)$ and $(ii)$, we compute the pullback to the fixed locus $Z$,
$$\iota_Z^{\star}\left(b_1^{j_1}b_1^{j_1+g} \cdot \ldots\cdot b_1^{j_s} b_1^{j_s+g}\right)=\prod_{i=1}^{s}
(y_1^{j_i}+\ldots+y_r^{j_i})(y_1^{j_i+g}+\ldots+y_r^{j_i+g})$$ $$=\prod_{i=1}^{s} \left(y_1^{j_i}
y_1^{j_i+g}+\ldots+y_r^{j_i}y_r^{j_i+g}\right)=\sum_{a_1+\ldots+a_r=s} \binom{s}{a_1, \ldots, a_r}
\frac{(g-a_1)!}{g!}\cdots \frac{(g-a_r)!}{g!} \theta_1^{a_1}\cdots \theta_r^{a_r}.$$ Here $a_i$ records the number of $y$s carrying the subscript $i$. 

Now following the exact same steps as in the proof of Lemma \ref{sumlemma}, we find $$\sum_{d_1+\ldots+d_r=d}
\int_{Z_{\{d_i\}}}\frac{\theta_1^{a_1} \cdot \ldots\cdot\theta_r^{a_r} \cdot R(\tilde {x_1}, \ldots,\tilde
{x_r})}{e_{T}(\n^{vir}_{Z})}=u\cdot \left(R J^{\gbar}\right)(\lambda_1, \ldots, \lambda_r)\cdot
\prod_{i=1}^{r}\frac{g!}{(g-a_i)!}\frac{\lambda_i^{a_i}}{N^{a_i}}$$ if $a_1+\ldots+a_r\leq d$ and $0$ otherwise. Proposition \ref{bees} follows by applying the binomial theorem to sum over the $a_i$s.\vskip.07in

It is harder to intersect $f$ classes. If one attempts the computation using degeneration methods, excess intersections
appearing on the boundary of the Quot scheme need to be considered. It does not seem possible to evaluate these intersections via the methods of \cite {bertram} \cite {sieberttian}. Here, we exemplify the power of the localization method by
calculating the intersection product \begin{equation}\label{newvi2} f_l \cdot P(a_1, \ldots, a_r) \cap
[Q_d]^{vir}.\end {equation} 

{\it Proof of Theorem $\ref {vinew}$.} By observation $(i)$ above, only the second term $$\sum_{q=1}^{r}d_q \sigma_{l-1;q}(\tilde x_1, \ldots, \tilde x_r)$$ in the formula for $f_l$, $2\leq l \leq r$ contributes to the 
evaluation
of $\eqref{newvi2}$. Let $$R(z_1, \ldots, z_r) = z_1^{\alpha_1} \cdots z_r^{\alpha_r}$$ be a
monomial of degree $e-l+1$.  As in Lemma \ref{sumlemma}, we compute \begin{equation} \sum_{d_1 + \cdots d_r = d} \int_{Z_{\{d_i\}}} \frac{R(\tilde 
x_1, \ldots, \tilde x_r) \left ( \sum_{q=1}^{r} d_q \sigma_{l-1;q}(\tilde x_1, \ldots \tilde x_r)  \right)  }{e_{T} (\n^{vir}_{Z_{\{d_i\}}})} = \end{equation} $$= u N^{r} \cdot R(\lambda_1, \ldots,
\lambda_r) J^{g-1} (\lambda_1, \ldots, \lambda_r) \cdot \sum_{d_1 + \cdots d_r = d}{\bar{I}}_{\{d_i\}, f_l},$$ where
$${\bar{I}}_{\{d_i\}, f_l} = \text {Res}_{\overline x_i=0} \left\{\prod_{1\leq i \leq r} \frac{(1+\bar{x}_i)^{\alpha_i +
(N-1)g}} {\left ( (1 + \bar{x}_i)^N -1 \right )^{d_i +1}} \prod_{1\leq i <j \leq r} \left ( 1 + \frac{\lambda_i \bar{x}_i
- \lambda_j \bar{x}_j}{\lambda_i - \lambda_j} \right )^{-2\gbar} \cdot\right . $$ $$\left. \cdot \sum_{q=1}^r d_q
\sigma_{l-1;q}\left(\lambda_1 (1 + \bar{x}_1), \ldots, \lambda_r (1 + \bar{x}_r)\right)\right \}.  $$ We evaluate the
coefficient of $$\prod_{i<j} (-1)^{k_{ij} - m_{ij}} \binom{-2\gbar}{k_{ij}}\binom{k_{ij}}{m_{ij}}
\frac{\lambda_{i}^{m_{ij}} \lambda_j^{k_{ij} - m_{ij}}}{(\lambda_i-\lambda_j)^{k_{ij}}}$$ in the sum $\sum_{d_1 + \cdots + d_r =
d} \bar{I}_{\{d_i \}, f_l}.$ We write $\lambda=(\lambda_1, \ldots, \lambda_r)$, and keep the same notation as in 
the proof of the Vafa-Intriligator formula. Summing over degrees just as in
Vafa-Intriligator, the coefficient in question equals $$\frac{1}{N^r}\prod_{1 \leq i \leq r}
\sum_{p_i = 0}^{m_i} (-1)^{m_i-p_i} \binom{m_i}{p_i} \binom{d-1 + \frac{s+n}{N}}{d-1}\left ( \sum_{q=1}^{r}
\frac{l_q + p_q}{N} \sigma_{l-1;q} ({\bf \lambda} ) \right ) = $$ $$ = \frac{1}{N^{r+1}}\sum_{n=0}^{k} (-1)^{k-n}
\binom{d-1 + \frac{s+n}{N}}{d-1} \left ( \sum_{q=1}^{r} \binom{k}{n} l_q \sigma_{l-1;q} ({\bf \lambda}) +
\sum_{q=1}^{r}\binom{k-1}{n-1} m_q \sigma_{l-1;q} ({\bf \lambda}) \right ).$$ By summing now over $m_{ij}$, then over $k_{ij}$, $n$ and $k$, we get
\begin{equation}
\label{sumfl}
\sum_{d_1 + \cdots + d_r = d} \bar{I}_{\{d_i \}, f_l} = \frac{1}{N^{r+1}} 
\left ( \sum_{q= 1}^{r}l_q \sigma_{l-1;q} ({\bf \lambda}) - 2\bar{g} \sum_{i<j} \frac{\lambda_i \sigma_{l-1;i}({\bf \lambda}) 
 - 
\lambda_j \sigma_{l-1;j}({\bf \lambda})}{\lambda_i - \lambda_j} \right ).
\end{equation}
The second sum in the right-hand side of \eqref{sumfl} simplifies to $$\frac{1}{2} (r-l)(r-l+1)\cdot 
\sigma_{l-1}.$$ 
Recalling the definition of $l_q = \alpha_q + (N-1) \gbar,$ we can finally write 
\begin{equation}
f_l \cdot P(a_1, \ldots, a_r) \cap [Q_d]^{vir} = 
\frac{u}{N} \sum_{\lambda_1, \ldots, \lambda_r} (\mathcal D_{l} R) (\lambda_1, \ldots, \lambda_r)\cdot J^{g-1}(\lambda_1, \ldots, \lambda_r).
\end{equation}
Here $\mathcal D_{l}$ denotes the differential operator given by the formula 
\begin {equation}\label{tl}\mathcal D_l R (z)=\gbar (r-l+1)(N-r+l-1)\cdot \sigma_{l-1} (z) \cdot R(z) + \sum_{q=1}^{r} \sigma_{l-1;q}(z) \cdot z_q \cdot \frac{\partial R}{\partial z_q}.
\end {equation}
Theorem $\ref {vinew}$ is a rewriting of the above equations. 

\begin{remark}We observe that a general product of $f$ classes in $\eqref{newvi2}$ yields differential operators of higher order.
Their exact expressions are not yet known to the authors, but we will return to this question in future work.
Nonetheless, let us note that certain intersection numbers involving higher powers of $f_2$ have been computed in \cite
{opreamarian} in rank $2$.\end{remark}

\section {An Application to the Moduli Space of Bundles} We would like to present another
application of the intersection theory of the $a$ classes on the Quot scheme. We prove Theorem
$\ref{modbun}$: we demonstrate the existence of a non-zero element of degree $r(r-1)\bar g$ in
the Pontrjagin ring of the moduli space $\mathcal M$ of rank $r$ bundles of degree $d$ on a genus
$g$ curve, where $\gcd(r, d)=1$. We will recover in this way Theorem $2$ in \cite {kirwan}.

To get started, let us write $\mathcal V$ for the universal bundle on $\mathcal M \times
C$. Let $p: \mathcal M \times C \to \mathcal M$ be the projection. We write $w_1, \ldots, w_r$ for the Chern roots of
the restricted bundle ${\mathcal V}_p$ on $\mathcal M\times \{p\}$ for some point $p \in C$. The
computation of the Pontrjagin class of $\mathcal M$ is well known \cite {N}, $$p(\mathcal
M)=\prod_{i<j}(1+(w_i-w_j)^2)^{2\gbar}.$$ In particular, the following element of algebraic
degree $r(r-1)\bar g$ $$\Theta = \prod_{i<j} (w_i-w_j)^{2\bar g}$$ is contained in the Pontrjagin
ring. We will assume for a contradiction that $\Theta = 0$.

We will move the computation from $\mathcal M$ to a suitable scheme $\quot$ making use of the
intermediate space  
$$\p_N(r, d)=\p(p_{\star} \mathcal V^{\oplus N})\stackrel{\pi} \longrightarrow\mathcal M.$$
This general setup is explained in \cite {marian} and \cite {opreamarian}; here, we only summarize the results we will need. The closed points of $\p_N(r, d)$ are pairs $(V, \phi)$ consisting of a stable
vector bundle of rank $r$ and degree $d$, and a nonzero morphism $$\phi: \mathcal O^N\to V.$$ Therefore,
$\p_N(r, d)$ and $\quot$ agree on the open subscheme where $V=E^{\vee}$ is stable and the morphism
$\phi$ is generically surjective. 

Just as for $\quot$, we can consider the K\"{u}nneth components of the universal bundles on $\p_N(r,
d)\times C$ and $\mathcal M\times C$, thus obtaining the corresponding $a$, $b$ and $f$ classes
on the moduli spaces $\p_N(r, d)$ and $\mathcal M$. The analogues of the class $\Theta$ can also
be defined on $\p_N (r,d)$ and $\quot.$

Tensoring with line bundles we may assume the degree $d$ is large (hence making $Q_d$ irreducible
and generically smooth). Moreover, we take $N$ sufficiently large such that the results of \cite
{marian} apply, and in particular the following equality of intersection numbers holds,
\begin{equation}\label{eqint}\int_{\quot} \Theta \cdot a_r^{M+\gbar} \cdot a_1 =
\int_{{\p_{N}(r, d)}} \Theta\cdot a_r^{M+\gbar} \cdot a_1=\int_{\mathcal M}\Theta \cdot \pi_{\star} (
a_r^{M+ \gbar} \cdot a_1)=0.\end{equation} 

We are free to pick any large $N$ we want, so we may assume $$Nd \equiv 1\mod r, \;\text { letting } rM=N(d-r\gbar)-1.$$ We will use the Vafa-Intriligator formula to see that the
left hand side of $\eqref{eqint}$ is non-zero. We need to compute the following sum over distinct
roots of unity $$\sum_{\lambda_1, \ldots, \lambda_r} \prod_{i<j} (\lambda_i-\lambda_j)^{2\gbar}
\cdot (\lambda_1 \ldots \lambda_r)^{M+\gbar} \cdot (\lambda_1+\ldots+\lambda_r)\cdot J^{g-1} (\lambda_1, \ldots, \lambda_r)=$$
$$=N^{r\gbar}\sum_{\lambda_1, \ldots, \lambda_r} (\lambda_1 \ldots \lambda_r)^{M} \cdot (\lambda_1+ \ldots
+ \lambda_r)=r N^{r\gbar} \sum_{\lambda_1, \ldots, \lambda_r} \lambda_1^{M+1} \lambda_2^{M}
\ldots \lambda_r^{M} =$$ $$=N^{r\gbar+1} \sum_{\mu_2, \ldots, \mu_{r}} \mu_2^{M} \ldots \mu_r^{M}.$$
Here, we set $\mu_i=\lambda_i\lambda_1^{-1}$. These are {\it distinct} $N$-roots of $1$ {\it not
equal to $1$}. Note that the constant factors change because of the repetitions introduced by this relabeling. 

Since $\mu\to \mu^{M}$ permutes the $N$-roots of unity, after writing $\zeta_i =
\mu_i^M$, we have to show $$\sum_{\zeta_2, \ldots, \zeta_r} \zeta_2 \cdot \ldots \cdot
\zeta_r\neq 0.$$ Letting $\zeta$ be a primitive root of unity, the expression above can be
evaluated as the coefficient of $t^{r-1}$ in the product $$(1+ \zeta t) \cdot \ldots
\cdot (1+ \zeta^{N-1} t) = 1+ t+ \ldots + t^{N-1}.$$ This gives the desired
contradiction.

\end{document}